\documentclass[11pt]{article}
\usepackage{}
\usepackage{amssymb}
\usepackage{anysize}
\marginsize{2.0cm}{2.0cm}{2.0 cm}{2.0 cm}
\usepackage{setspace}
\usepackage{amsmath}
\usepackage{amsthm}
\usepackage{amsfonts}
\usepackage{mathrsfs}
\newtheorem{thm}{Theorem}[section]
\newtheorem{prop}{Proposition}[section]
\newtheorem{lem}{Lemma}[section]
\newtheorem{rem}{Remark}[section]
\newtheorem{cor}{Corollary}[section]

\newtheorem{exm}{Example}[section]

\newtheorem{con}{Condition}[section]
\newtheorem{example}{Example}[section]

\numberwithin{equation}{section} \allowdisplaybreaks[4]
\begin{document}

\title{Finite horizon risk-sensitive continuous-time Markov decision processes with unbounded transition and cost rates}

\author{Xin Guo\thanks{Department of Mathematical Sciences, University of
Liverpool, Liverpool, L69 7ZL, U.K.. E-mail:  X.Guo21@liv.ac.uk.},~ Qiuli Liu\thanks{School of Mathematical Sciences, South China Normal University, Guangzhou, 510631. China. E-mail: liuql2007@aliyun.com.}~ and~ Yi
Zhang\thanks{Department of Mathematical Sciences, University of
Liverpool, Liverpool, L69 7ZL, U.K.. E-mail: yi.zhang@liv.ac.uk.}}

\date{}
\maketitle
\par\noindent{\bf Abstract:}  We consider a risk-sensitive continuous-time Markov decision process over a finite time duration. Under the conditions that can be satisfied by unbounded transition and cost rates, we show the existence of an optimal policy, and the existence and uniqueness of the solution to the optimality equation out of a class of possibly unbounded functions, to which the Feynman-Kac formula was also justified to hold.
\bigskip

\bigskip


\section{Introduction} \label{intro}

Risk-sensitive Markov decision processes (in discrete-time) have been studied intensively since 1970s, with one of the pioneering works being \cite{Howard:1972}, and a recent and updated work being \cite{BauerleRieder:2014}, to which the interested reader is referred for more references. Compared to the discrete-time framework, there have been fewer works on risk-sensitive CTMDPs (continuous-time Markov decision processes), also known as controlled Markov pure jump processes. An early work on this topic seems to be \cite{Piunovski:1985}, which obtained verification theorems and solved in closed-form meaningful examples of problems over a fixed time duration. In the recent years, there have been reviving interests in risk-sensitive CTMDPs, see \cite{GS-14,GuoZhangAMO:2017,Wei:2016} for problems with a finite horizon, \cite{Zhang:2017} for problems over an infinite horizon,  \cite{GS-14,KC-13,WeiC:2016} for problems with average criteria, and \cite{BauerleP:2017} for an optimal stopping problem with a more general utility function than the exponential one.

In greater detail, the CTMDP considered in \cite{GS-14} is with bounded transition and cost rates. In \cite{Wei:2016}, the boundedness on the transition rate was relaxed and replaced by a drift-type condition, but the cost rate was still assumed to be bounded. Both papers followed the same line of reasoning: they showed the existence of a solution to the optimality equation, and then showed that the solution coincides with the value function of the problem by applying the Feynman-Kac formula. In Section 7 of \cite{Wei:2016}, the author mentioned that following his method it was unclear how to relax the boundedness assumption on the cost rate at that time, as a suitable version of the Feynman-Kac formula must be established first. The present paper provides a response to this. In greater detail, the main contributions are the following. We provide conditions that allow unbounded transition and (not necessarily nonnegative) cost rates, under which a suitable version of the Feynman-Kac formula was established, and we show that the value function is the unique solution out of a large enough class of functions (possibly unbounded with unbounded derivatives with respect to time) to the optimality equation. It is important for practical applications to consider models with unbounded transition and cost rates. We illustrate this with an example of controlled $M/M/\infty$ queueing system. Compared with \cite{GS-14,Wei:2016}, which concentrated on Markov policies, we consider a more general class of policies. When the cost rate is nonnegative,  a different method was followed in \cite{GuoZhangAMO:2017}, which is not based on the Feynman-Kac formula. If the cost rate is nonnegative, then the conditions on the transition and cost rates in \cite{GuoZhangAMO:2017} are weaker than in the present paper. Moreover, in that general setup of \cite{GuoZhangAMO:2017}, the value function is generally not the unique solution to the optimality equation. In this sense, the present paper also complements \cite{GuoZhangAMO:2017}.

The rest of the paper is organized as follows. In Section \ref{section} we describe the optimal control problem under consideration. Section \ref{sectionPre} contains preliminary results, where we establish a version of the Feynman-Kac formula. The optimality results are proved in Section \ref{sectionMS}. This paper is finished with a conclusion in Section \ref{sectionConclusion}.

\section{Model description}\label{section}

\textbf{Notation:} For a Borel space $X$ endowed with the Borel
$\sigma$-algebra ${\cal B}(X)$, we denote by $\mathbb{C}_b(X)$ the space of all bounded continuous functions on $X$ with the norm $\|u\|:=\sup_{x\in X}|u(x)|$. Throughout this paper, measurability is understood in the Borel sense.

We consider the CTMDP model
$
{\cal M}:=\{S, A, A(\cdot,\cdot), q, c, g\}
$
consisting of the following elements.
The state space $S$ is a denumerable set, endowed with the discrete topology. The action space $A$ is a (nonempty) Borel space. The multifunction $(t,i)\in [0,\infty)\times S\rightarrow A(t,i)\in {\cal B}(A)$ specifies the set of admissible action spaces given the current time and state, and is assumed to be with a measurable graph $K:=\{(t,i,a)\in [0,\infty)\times S\times A: a\in A(t,i)\}$, containing the graph of some measurable mapping from $[0,\infty)\times S$ to $A$.
The transition rate is given by a signed
kernel $q$ on $S$ given $K$, assumed to satisfy $q(j|t,i,a)\ge 0$ if $j\not=i$ with $j,i\in S,$
$q(S|t,i,a)\equiv 0,$ and
 \begin{eqnarray}\label{Q}
q^{*}(i):=\sup_{t\geq 0, a\in A(t,i)}q(t,i,a)<\infty,~\forall~  i\in  S,
\end{eqnarray}
where $q(t,i,a):=-q(i|t,i,a)\geq0$ for all $(t,i,a)\in K$. The running cost rate $c$ is a measurable function on $K$. We shall consider the problem over a finite time duration. The terminal cost $g$ is a function on $S$.

We briefly describe the construction of a CTMDP as in \cite{Kitaev:1986,Kit95}. Let
$S_\Delta:=S\bigcup\{\Delta\}$ (with some $\Delta\not\in S$ being an isolated point),
$\Omega^{0}:=(S\times(0,\infty))^{\infty}$ be the countable product. The canonical sample space
$\Omega$ is the union of $\Omega^0$ and all the sequences in the form of $(i_{0},\theta_{1},i_{1},\ldots,\theta_{k},i_{k},
\infty,\Delta,\infty, \ldots)$ for some $k\ge 0$ (accepting $\theta_0:=0$).  Let
$\mathcal{F}$ be the Borel $\sigma$-algebra on $\Omega$. For each $\omega\in \Omega$, introduce
$T_{0}(\omega):=0,$  $T_{k+1}(\omega):=\theta_{1}+\theta_{2}+\ldots+\theta_{k+1},$  $X_k(\omega):=i_k.$
In what follows, the argument $\omega$ is often omitted. Let ${\cal F}_t$ be the internal history of the marked point process $\{T_n,X_n\}$. Let
$T_{\infty}:=\lim_{k\rightarrow\infty}T_{k}$. The controlled
process $\{\xi_t\}$ is defined by
\begin{eqnarray*}
 \xi_t(\omega):=\sum_{k\geq0}I_{\{T_{k}\leq t<T_{k+1}\}}i_{k}+\Delta I_{\{t\geq T_\infty\}}, \forall~\ t\geq 0.
\end{eqnarray*}
Here and below, $I_{E}$ stands for the indicator function on any
set $E$, and for notational convenience, we defined that $i\cdot 0=0$ and $i\cdot 1=i$ for each $i\in S_\Delta.$

We do not intend to consider the
controlled process after moment $T_\infty$, and put
\begin{eqnarray*}
q(\cdot|t,\Delta,a_\Delta):\equiv 0,~ r(t,\Delta,a_\Delta):\equiv
0,~ A(t,\Delta):=\{a_\Delta\},~ A_\Delta:=A\cup\{a_\Delta\},
\end{eqnarray*}
where
$a_\Delta\notin A$ is an isolated point.

A  (history-dependent) policy $\pi$  is determined and often identified by a sequence of stochastic kernels $\{\pi^k,k\geq 0\}$ such that
\begin{eqnarray*}
\pi(da|\omega,t)&=&I_{\{t=0\}}\pi^{0}(da|i_{0},0)+\sum_{k\geq0}I_{\{T_{k}<t\leq T_{k+1} \}}\pi^{k}(da|i_{0},\theta_{1},i_{1},\ldots,\theta_{k},i_{k},t-T_{k}) \nonumber \\
&&   + I_{\{t\geq T_{\infty}\}}I_{\{a_\Delta\}}(da).
\end{eqnarray*}
A policy $\pi$ is called Markov if, with slight abuse of notations, $\pi(
da |\omega,t)=\pi( da
|\xi_{t-},t)$, which is denoted by $\pi_{t}(da|\cdot)$, where $\xi_{t-}=\lim_{s\uparrow t} \xi_s.$ A Markov policy
$\pi_{t}(da|\cdot)$ is called deterministic Markov if
there exists a measurable mapping
$f$ on $[0,\infty)\times S$ such that $\pi_{t}(da|i)$ is a
Dirac measure concentrated at $f(t,i)$. A deterministic Markov policy will
be denoted by the underlying measurable mapping $f$.
We denote by $\Pi$ the set of all
policies, by $\Pi_{m}^r$ the set of all Markov policies, and by $\Pi_{m}^d$ the set of all deterministic Markov policies.

For each $\pi\in \Pi,$ the random measure $m^\pi$ defined by
\begin{eqnarray}
m^\pi(j|\omega,t)dt:=\int_{A}q(j\setminus \{\xi_{t-}\}|t,\xi_{t-},a)\pi(da|\omega,t)dt \label{RM}
\end{eqnarray}
is predictable, see \cite{Jac75}. For each $\pi\in \Pi$ and $i\in S$, let $\mathbb{P}_i^\pi$ be the probability on $(\Omega,{\cal F})$ such that $\mathbb{P}_{i}^{\pi}(\xi_{0}=i)=1$, and with respect to which, $m^\pi(j|\omega,t)dt$ is the dual predictable projection of the random measure $\sum_{n\ge 1}\delta_{(T_n,X_{n})}(dt,dx)$ of the marked point process $\{T_n,X_n\}$ on ${\cal B}((0,\infty)\times S)$, see \cite{Jac75,Kitaev:1986} or Chapter 4 of \cite{Kit95} for more details. Let $\mathbb{E}_i^\pi$ be the expectation taken with respect to $\mathbb{P}_i^\pi$.

For the intuitive description, a CTMDP is a continuous-time Markov pure jump process whose local characteristics (transition intensity and post-jump distributions) are controlled. After the $n$-th jump, and a history of state and sojourn times $h_n=(i_0,\theta_1,\dots,\theta_n,i_n)$ is observed with $\theta_n<\infty$, the conditional (joint) distribution of the next state and sojourn time is determined by
$\int_A q(j|t_n+t,i_n,a)\pi_n(da|h_n,t-t_n)e^{-\int_0^t \int_A q(s+t_n,i_n,a)\pi_n(da|h_n,s)ds}dt$, $j\ne i_n,$ where $t_n$ is the observed value of the $n$-th jump moment. In particular, the next sojourn time has the conditional distribution obeying a nonstationary exponential distribution, and a policy specifies the selection of an action at any time moment based on the observed history.

We consider the following optimal control problem over the finite time duration $T>0$:
\begin{eqnarray}\label{New02}
\mbox{Minimize over $\pi\in\Pi$:~} {\cal V}(\pi,i):=\mathbb{E}_{i}^{\pi}\left[e^{\int_{0}^{T}\int_{A}c(t,\xi_{t},a)\pi(da|\omega,t)dt+g(\xi_T)}\right].
\end{eqnarray}
Conditions imposed in the next section guarantee that the above expectation and integral are well defined.
For each $i\in S$, let
\begin{eqnarray*}
{\cal V}^*(i)=\inf_{\pi\in\Pi}{\cal V}(\pi,i).
\end{eqnarray*}
A policy $\pi^{*}\in\Pi$ is said to be optimal if ${\cal V}(\pi^*,i)={\cal V}^*(i)$ for all $i\in S$.

Problem (\ref{New02}) is often said to be with a risk-sensitive criterion, as the exponential utility reflects that the decision maker is increasingly averse to the higher cost, see \cite{Howard:1972}. This is in contrast with a linear utility, which is called risk-neutral. In discrete-time, risk-sensitive Markov decision processes received increasing interest in the recent years, see \cite{Cavazos-Cadena:2000a,Cavazos-Cadena:2000b,Jaskiewicz:2008,Patek:2001} for example. These works mainly consider infinite-horizon problems; in the discrete-time setup, problems on finite horizon can be readily solved using backward induction. See also \cite{BauerleRieder:2014}, which considered a more general utility function.

The objective of this paper is to provide conditions that can be satisfied by unbounded transition and cost rates, under which, there exists a deterministic Markov optimal policy, and the optimality equation has a unique solution out of a certain class of functions. We present an example in the next section, demonstrating a natural application of CTMDPs to controlled queueing system, where the transition and cost rates are both unbounded and thus not covered by the previous literature.

\section{Preliminaries}\label{sectionPre}

In this section, we impose a set of conditions allowing one to consider unbounded transition and cost rates, see Example \ref{Example01} below, and present several preliminary statements, which will serve the proof of Theorem \ref{T-4.1} below.
\begin{con}\label{ass:3.1}
There exist a $[1,\infty)$-valued function $V$ defined on $S$ and constants $\rho>0$, $M>1$ such that
\begin{itemize}
\item[(a)]  $\sum_{j\in S}q(j|t,i,a)V(j)\leq \rho V(i)$  for each $(t,i,a)\in  K$;

\item[(b)] $q^*(i)\leq MV(i)$ for all $i\in S$, where $q^*(i)$ is as in (\ref{Q});

\item[(c)]  $e^{2(1+T)|c(t,i,a)|}\leq MV(i)$ for each $(t,i,a)\in K$, and $e^{2(1+T)|g(i)|}\leq MV(i)$ for each $i\in S$.
\end{itemize}
\end{con}
The immediate and relevant consequences of Condition \ref{ass:3.1} are collected in the next lemma.
\begin{lem}\label{lem3.1} Suppose Condition  \ref{ass:3.1} is satisfied. For each $\pi\in\Pi$, the following assertions hold.
\begin{itemize}
\item[(a)]  $\mathbb{P}_{i}^{\pi}(T_\infty=\infty)=1$ for each $i\in S$.
\item[(b)] $\mathbb{E}_{i}^{\pi}[V(\xi_t)] \leq e^{\rho t}V(i)$, for each $t\geq 0$ and $i\in S$.
\item[(c)]  $  V(\pi,i)\leq Me^{T\rho}V(i) $ for all $i\in S$ and $\pi\in\Pi$.
\end{itemize}
\end{lem}

\par\noindent\textit{Proof.} Parts (a) and (b) are known, see e.g., \cite{GuoABP:2010,PiunovskiZ:2011,ABPZY:2014}. We next verify part (c). By part (a), for $\mathbb{P}_i^\pi$-almost all $\omega\in\Omega,$ there are finitely many values taken by in $\{\xi_t(\omega)\}$ over $[0,T]$. For such $\omega\in \Omega,$ by Condition \ref{ass:3.1}(c), we legitimately write
\begin{eqnarray*}
\int_0^T \int_A  c(t,\xi_t,a)\pi(da|\omega,t)dt+g(\xi_T)=\int_{(0,T]} \int_A \tilde c(t,\xi_t,a)\pi(da|\omega,t)\mu(dt),
\end{eqnarray*}
where $\mu(dt)= I_{[0,T)}(t)dt + \delta_{{T}}(dt)$, with $\delta_T(dt)$ being the Dirac measure concentrated on $\{T\}$, and $\tilde c(t,i,a):= c(t,i,a)I_{[0,T)}(t)+g(i)I_{\{T\}}(t)$ for each $(t,i,a)\in K.$  Now,
\begin{eqnarray}\label{H}
&&\mathbb{E}_{i}^{\pi}\left[e^{\int_0^T\int_A  c(t,\xi_t,a)\pi(da|\omega,t)dt+g(\xi_T)}\right]=\mathbb{E}_{i}^{\pi}\left[e^{\int_{[0,T]}\int_A (1+T)\tilde c(t,\xi_t,a)\pi(da|\omega,t)\frac{\mu(dt)}{T+1}}\right] \nonumber\\
&\le& \mathbb{E}_{i}^{\pi}\left[\frac{1}{1+T}\int_{[0,T]} e^{(1+T)\int_A |\tilde c(t,\xi_t,a)|\pi(da|\omega,t)}\mu(dt)\right] \leq \frac{M}{1+T}\mathbb{E}_{i}^{\pi}\left[\int_0^T V(\xi_t)dt+V(\xi_T) \right] \nonumber \\
&\leq&  Me^{\rho T}V(i),
\end{eqnarray}
where the first inequality is by the Jensen inequality, the second inequality is by Condition \ref{ass:3.1}(c), and the last inequality is by part (b). $\hfill\Box$
\bigskip

Part (a) of the previous lemma asserts that under the imposed conditions therein, the controlled process is nonexplosive under each policy. This fact is used in the proof of Theorem \ref{Th-3.2} below, see the first paragraph therein as well as (\ref{New01}).
\begin{con}\label{ass:3.1b}
 There exist a $[1,\infty)$-valued function  $V_1$ defined on $S$, and constants
$\rho_1>0$, $M_1>0$ such that
\begin{itemize}
\item[(a)] $\sum_{j\in S}V_1^2(j)q(j|t,i,a)\leq \rho_1V_1^2(i)$ for each $(t,i,a)\in  K$;

\item[(b)] $V^2(i)\leq M_1V_1(i)$ for all $i\in S$, with the function $V$  as the Condition \ref{ass:3.1}.
\end{itemize}
\end{con}
The role of this condition is seen in the proof of Theorem \ref{Th-3.2}, where the Cauchy-Schwarz inequality is used, see (\ref{GDA}) therein. Conditions \ref{ass:3.1} and \ref{ass:3.1b} guarantee the growth of the value function and its derivative to be suitably bounded by the weight functions $V$ and $V_1$, and it is out of this class of functions that we show the Feynman-Kac formula applies. The previous works \cite{GS-14,Wei:2016} only showed that the Feynman-Kac formula is applicable to a class of bounded functions, and so confined themselves to the class of bounded cost rates, which excludes some potentially interesting applications. Let us formulate such an example, which are with unbounded transition and cost rates and satisfy Conditions \ref{ass:3.1} and \ref{ass:3.1b}.
\begin{example}\label{Example01}
Consider a controlled $M/M/\infty$ queueing system, where the common service rate $a$ of each server can be tuned from a finite interval $[\underline{\mu},\overline{\mu}]\subseteq [0,\infty]$. Let the arrival rate be denoted by $\lambda>0$. The holding cost is $C_1 i$ given the current number of jobs in the system being $i\ge 0$, where $C_1>0$ is a constant, and maintaining a service rate at $\mu$ costs $\mu$ per unit time. A terminal reward of $C_2 i$ is received if there are $i$ jobs remaining in the system at the end of the horizon $[0,T]$, where $C_2\in(-\infty,\infty)$ is a constant. The decision maker aims at the optimal control of the service rate to minimize the expected exponential utility of the total cost over the horizon $[0,T]$.

This problem can be formulated as a CTMDP with the following primitives. The state space is $S=\{0,1,\dots\}$, the action space is $[\underline{\mu},\overline{\mu}]\equiv A(t,i)$. The transition rate is given by $q(i+1|t,i,a)\equiv \lambda$, $q(i-1|t,i,a)=a i$ if $i\ge 1$, $q(t,i,a)=\lambda+ a i$ if $i>0$, and $q(t,0,a)=\lambda.$ The running cost rate is given by $c(t,i,a)= C_1i+a$, and the terminal cost is given by $g(i)=-C_2 i.$

Observe the following. Let $d>0$ be a fixed constant. Let $\rho(d): =e^{d+1}\lambda.$ Then for each constant $\rho\ge \rho(d),$ $\sum_{j\in S}q(j|t,i,a)e^{d j}=e^{d (i+1)}\lambda +e^{d(i-1)}a-(\lambda+a)e^{di}\le \rho e^{d i}$ for each $i\ge 1$, and $\sum_{j\in S}q(j|t,0,a)e^{d j}=\lambda e^{d}-\lambda \le \rho $. Therefore, for the verification of Condition \ref{ass:3.1}, one can take
$M= e^{2(1+T)\overline{\mu}}+\overline{\mu}+\lambda$, $V(i)=e^{d_1 i}$ with $d_1=2(1+T)(C_1+|C_2|)$, $\rho=\rho(d_1)$. For the verification of Condition \ref{ass:3.1b}, one can take $M_1=1$, and $V_1(i)=e^{d_2 i}$ with $d_2= 2 d_1$, and $\rho_1=\rho(d_2)$.
\end{example}

Let us introduce some additional notations, which will be needed in the next statement. In particular, it formalizes what we mean in the Introduction by ``a large enough class of functions'' to which, the Feynman-Kac formula applies.
Let $X$ be a Borel space. For each measurable function $\psi$ on $[0,T]\times X$, if $\psi(\cdot,x)$ is absolutely continuous on $[0,T]$, then we put $\psi'$ a measurable function on $[0,T]\times X$ such that $\psi(t,x)-\psi(0,x)=\int_{0}^t \psi'(s,x)ds$ for each $x\in X$ and $t\in[0,T]$. Consider the functions $V$ and $V_1$ as in Conditions \ref{ass:3.1} and \ref{ass:3.1b}. A function $\varphi$ on $[0,T]\times S$ is called $V$-bounded if the $V$-weighted norm of $\varphi$, $\|\varphi\|_{V}:=\sup_{(t,i)\in [0,T]\times S}\frac{|\varphi(t,i)|}{V(i)}$, is finite. Let
$C_{V,V_1}^{1}([0,T]\times S)$ be the collection of $V$-bounded functions $\varphi$ on $[0,T]\times S$ such that $\varphi(\cdot,i)$ is absolutely continuous on $[0,T]$ for each $i\in S,$ which admits some $V_1$-bounded $\varphi'$.

\begin{thm}\label{Th-3.2}
 Suppose Conditions \ref{ass:3.1} and \ref{ass:3.1b} are satisfied. Then, for each $i\in S$, $\pi\in\Pi$ and $\varphi\in C_{V,V_1}^{1}([0,T]\times S)$,
 \begin{eqnarray*}
 && \mathbb{E}_{i}^{\pi}\left[\int_{0}^{T}
\left(\psi'(\omega,t,\xi_{t})+\sum_{j\in S}\psi(\omega,t,j)\int_{A} q(j|t,\xi_t,a)\pi(da|\omega,t)\right)dt
\right] \\
&=&\mathbb{E}_{i}^{\pi}\left[\psi(\omega,T,\xi_{T})\right]-\varphi(0,i),
\end{eqnarray*}
where outside a $\mathbb{P}_i^\pi$-null set, say $\Omega\setminus \Omega',$ $T_\infty=\infty,$
\begin{eqnarray*}
\psi(\omega,t,j)=e^{\int_0^t \int_A c(v,\xi_v,a)\pi(da|\omega,v)dv}\varphi(t,j),~\forall~t\in[0,T],~j\in S,
\end{eqnarray*}
$\psi(\omega,\cdot,j)$ is absolutely continuous on $[0,T]$ so that we can take
\begin{eqnarray}\label{G-11}
    \psi'(\omega,t,j)&=&\int_A c(t,\xi_t,a)\pi(da|\omega,t)e^{\int_0^t\int_A c(v,\xi_v,a)\pi(da|\omega, v) dv}\varphi(t,j)\nonumber\\
     &&+e^{\int_0^t \int_A c(v,\xi_v,a)\pi(da|\omega,v)dv}\varphi'(t,j),
  \end{eqnarray}
for each $\omega\in \Omega'$ and $j\in S$.
\end{thm}
\par\noindent\textit{Proof.} According to Lemma \ref{lem3.1}(a), we concentrate on $\Omega'$ on which $T_\infty=\infty,$ and hence (\ref{G-11}) holds. Since $\varphi\in C_{V,V_1}^{1}([0,T]\times S)$, we have $|\varphi(t,i)|\leq \|\varphi\|_{V}V(i)$ for all $(t,i)\in [0,T]\times S$, which, together with the relation $(1+T)|c(v,i,a)|\leq MV(i)$ (by Condition \ref{ass:3.1}(c)), leads to
\begin{eqnarray}\label{G-1}
&&\left|\psi'(\omega,t,\xi_{t})\right|\nonumber\\
&\leq & \frac{M}{1+T}V(\xi_t)e^{\int_0^t \int_A |c(v,\xi_v,a)| \pi(da|\omega,v)dv}\|\varphi\|_{V}V(\xi_t) +\|\varphi'\|_{V_1}e^{\int_0^t\int_A|c(v,\xi_v,a)|\pi(da|\omega,v)dv}V_1(\xi_t),  \nonumber \\
&\leq& \frac{\|\varphi\|_{V}+\|\varphi'\|_{V_1}}{1+T}(1+T+MM_1)e^{\int_0^t\int_A|c(v,\xi_v,a)|\pi(da|\omega,v)dv}V_1(\xi_t).
\end{eqnarray}
By the Cauchy-Schwarz inequality,
\begin{eqnarray}\label{GDA}
       && \mathbb{E}_{i}^{\pi}\left[ e^{\int_0^t\int_A |c(v,\xi_v,a)|\pi(da|\omega,v)dv}V_1(\xi_t)\right]\leq  \sqrt{\mathbb{E}_{i}^{\pi}\left[e^{2\int_0^t\int_A|c(v,\xi_v,a)|\pi(da|\omega,v)dv}\right]\mathbb{E}_{i}^{\pi}\left[V_1^2(\xi_t)\right]} \nonumber \\
   &\leq & \mathbb{E}_{i}^{\pi}\left[e^{2\int_0^t\int_A|c(v,\xi_v,a)|\pi(da|\omega,v)dv}\right]\mathbb{E}_{i}^{\pi}\left[V_1^2(\xi_t)\right]\le  Me^{T\rho}V(i)\mathbb{E}_{i}^{\pi}\left[V_1^2(\xi_t)\right]\nonumber\\
   &\le&Me^{T\rho}V(i) e^{\rho_1 T}V_1^2(i) ,~t\in[0,T],
  \end{eqnarray}
where the second to the last inequality is obtained by a similar argument to the one for (\ref{H}), and the last inequality is by Lemma \ref{lem3.1}(b).
Now it follows from (\ref{G-1}) that
\begin{eqnarray}
\mathbb{E}_{i}^{\pi} \left[\int_0^{T}|\psi'(\omega,t,\xi_t)|dt\right]<\infty.
\label{F2}
\end{eqnarray}

 On the other hand,  by Conditions \ref{ass:3.1}  and \ref{ass:3.1b}, we have
   \begin{eqnarray*}
     &&\sum_{j\in S}e^{\int_0^t\int_A|c(v,\xi_v,a)|\pi(da|\omega,v)dv}|\varphi(t,j)|\left|\int_Aq(j|t,\xi_t,a)\pi(da|\omega,t)\right|\\
      &\leq &\|\varphi\|_{V}  \left(\rho V(\xi_t)+2MV^2(\xi_t)\right) e^{\int_0^t\int_A |c(v,\xi_v,a)|\pi(da|\omega,v)dv} \nonumber \\
     &\leq & \|\varphi\|_{V}M_1(\rho+2M)e^{\int_0^t\int_A | c(v,\xi_v,a)|\pi(da|\omega,v)dv}V_1(\xi_t).
   \end{eqnarray*}
Now it follows from (\ref{GDA}) that
\begin{eqnarray}\label{F1}
 \int_0^{T}\sum_{j\in S}\mathbb{E}_{i}^{\pi}\left[\left|\int_{A}q(j|t,\xi_t,a)\pi(da|\omega,t)\right||\psi(\omega,t,j)|\right]dt <\infty.
\end{eqnarray}

For each $0\le s\le T$,
\begin{eqnarray}\label{New01}
\psi(\omega,T,\xi_T)=\psi(\omega,0,\xi_0)+\int_{0}^T \psi'(\omega,t,\xi_t)dt+\sum_{n\ge 1}\int_{(0,T]}\Delta \psi(\omega,t,\xi_t)\delta_{T_n }(dt)
\end{eqnarray}
with $\Delta \psi(\omega,t,\xi_t):=\psi(\omega,t,\xi_t)-\psi(\omega,t-,\xi_{t-}).$ (Recall that the function $\psi(\omega,t,j)$ is absolutely continuous in $t$ over finite interval, and for each fixed $\omega\in\Omega'$ with $\Omega'$ being defined in the beginning of this proof, $\xi_t(\omega)$ is piecewise constant in $t\in[0,T]$, and assumes finitely many values over that interval.)
By (\ref{F2}) and (\ref{F1}), we take legitimately the expectation on the both sides of the previous equality, and obtain
\begin{eqnarray*}
&& \mathbb{E}_i^{\pi}\left[ \psi(\omega,T,\xi_T) \right]=\mathbb{E}_i^{\pi}\left[\psi(\omega,0,\xi_0)\right]+\mathbb{E}_{i}^{\pi}\left[ \int_{0}^T \psi'(\omega,t,\xi_t)dt \right]\\
&& + \mathbb{E}_{i}^{\pi}\left[ \sum_{n\ge 1}\int_{(0,T]}\Delta \psi(\omega,t,\xi_t)\delta_{T_n }(dt)\right] \\
 &=&\varphi(0,i)+\mathbb{E}_{i}^{\pi}\left[\int_{0}^T \psi'(\omega,t,\xi_t)dt\right]\\
 &&+\mathbb{E}_{i}^{\pi}\left[ \sum_{j\in S}\int_{(0,T]}(\psi(\omega,t,j)-\psi(\omega,t,\xi_{t-})) m^{\pi}(j|\omega,t)dt\right] \nonumber\\
&=&\varphi(0,i)+\mathbb{E}_{i}^{\pi}\left[ \int_{0}^T \psi'(\omega,t,\xi_t)dt\right]\nonumber\\
&&+\mathbb{E}_{i}^{\pi}\left[ \sum_{j\in S} \int_{0}^T\int_A \psi(\omega,t,j)q(j|t,\xi_{t-},a) \pi(da|\omega,t)dt\right],
\end{eqnarray*}
where the last equality holds because the random measure $m^\pi$ defined by (\ref{RM}) is the dual predictable projection of the random measure $\sum_{n\ge 1}\delta_{(T_n,X_{n})}(dt,dx)$ on ${\cal B}((0,\infty)\times S)$  under $\mathbb{P}_i^\pi$, see p.131 of \cite{Kit95}. The statement is proved. $\hfill\Box$
\bigskip

The above Feynman-Kac formula in the above theorem was justified in \cite{Wei:2016}, see Theorem 3.1 therein, when $\pi$ is a Markov policy, and $\varphi$ is assumed to be bounded.

The next statement provides a verification theorem, which was known in \cite{Piunovski:1985} when the transition rate is bounded.
\begin{cor}\label{ThC-3.2}
Suppose Conditions \ref{ass:3.1} and \ref{ass:3.1b} are satisfied. If there  exists   $\varphi\in C_{V,V_1}^{1}([0,T]\times S)$ and a deterministic Markov policy $f\in \Pi_m^d$ such that
\begin{eqnarray}\label{E-1}
\varphi(s,i)-e^{g(i)}&=&\int_s^T \inf_{a\in A(t,i)}\left\{c(t,i,a)\varphi(t,i)+\sum_{j\in S}\varphi(t,j)q(j|t,i,a)\right\} dt \nonumber\\
&=&\int_s^T  \left\{c(t,i,f(t,i))\varphi(t,i)+\sum_{j\in S}\varphi(t,j)q(j|t,i,f(t,i))\right\}dt,\nonumber\\
&&~s\in[0,T],~i\in S,
\end{eqnarray}
then
\begin{eqnarray}\label{GLZ1}
{\cal V}(f,i)= \varphi(0,i)={\cal V}^\ast(i),~\forall~i\in S.
\end{eqnarray}
\end{cor}
\par\noindent\textit{Proof.} Concentrate on $\Omega'$ as in the proof of the previous theorem. It holds for almost all $t\in[0,T]$ that
\begin{eqnarray*}
&&0=\varphi'(t,\xi_t)+\inf_{a\in A(t,\xi_t)}\left\{c(t,\xi_t,a)\varphi(t,\xi_t)+\sum_{j\in S}\varphi(t,j)q(j|t,\xi_t,a)\right\}\\
&=&\varphi'(t,\xi_t)+c(t,\xi_t,f(t,\xi_t))\varphi(t,\xi_t)+\sum_{j\in S}\varphi(t,j)q(j|t,\xi_t,f(t,\xi_t))\\
&\le&\varphi'(t,\xi_t)+\int_A \left\{c(t,\xi_t,a)\varphi(t,\xi_t)+\sum_{j\in S}\varphi(t,j)q(j|t,\xi_t,a)\right\}\pi(da|\omega,t).
\end{eqnarray*}
Now by applying Theorem \ref{Th-3.2} to the deterministic Markov policy $f$ and an arbitrarily fixed $\pi\in \Pi$, we see
\begin{eqnarray*}
&&{\cal V}(\pi,i)-\varphi(0,i)=\mathbb{E}_i^\pi\left[e^{\int_0^T \int_A c(v,\xi_v,a)\pi(da|\omega,v)dv}\varphi(T,\xi_T)  \right]-\varphi(0,i)\\
&=&\mathbb{E}_i^\pi\left[\int_0^T e^{\int_0^t \int_{A} c(v,\xi_v,a)\pi(da|\omega,v)dv}\int_A (c(t,\xi_t,a)\varphi(t,\xi_t)+\varphi'(t,\xi_t)+\sum_{j\in S}\varphi(t,j)q(j|t,\xi_t,a))\pi(da|\omega,t)\right]\\
&\ge&0,
\end{eqnarray*}
where the first equality holds because $\varphi(T,i)=e^{g(i)}$, see (\ref{E-1}); similarly, replacing $f$ for $\pi$ in the equalities in the above,
${\cal V}(f,i)-\varphi(0,i)=0.$ Consequently,
${\cal V}(f,i)=\varphi(0,i)\le {\cal V}(\pi,i)$ for each $i\in S.$ Since $\pi$ was arbitrarily fixed, ${\cal V}(f,i)=\varphi(0,i)={\cal V}^\ast(i),$ as required.
 $\hfill\Box$
\bigskip

According to the previous statement, (\ref{E-1}) is called the optimality equation, and the policy $f$ in (\ref{GLZ1}) is optimal.

The next statement was basically obtained in Theorem 2.1 in \cite{GS-14}, see also \cite{Wei:2016}.
\begin{prop}\label{ThCC-3.2}
Suppose that the transition and cost rates are bounded, i.e.,
\begin{eqnarray*}
\sup_{i\in S}q^*(i)<\infty,~\sup_{(t,i,a)\in K}|c(t,i,a)|<\infty,~\sup_{i\in S}|g(i)|<\infty.
\end{eqnarray*}
If for each $i\in S$ and $t\in[0,T],$ $A(t,i)$ is compact, $c(t,i,a)$ is lower semicontinuous in $a\in A(t,i)$, and $q(j|t,i,a)$ is continuous in $a\in A(t,i)$,
then there exists a unique $\varphi$ in $C^{1}_{1,1}([0,T]\times S)$ and some $f\in \Pi_m^d$ satisfying (\ref{E-1}) and (\ref{GLZ1}).
\end{prop}
The main objective in this paper is to relax the boundedness requirements in the previous statement.

\section{Optimality result}\label{sectionMS}
We impose the following condition, which guarantees the existence of an optimal policy.
 \begin{con}\label{ass:3.3}
\begin{itemize}
\item[(a)] For each $(t,i)\in [0,T]\times S$, $A(t,i)$ is compact.
\item[(b)] For each $t\in [0,T], i,j\in S$, the function  $q(j|t,i,a)$  is continuous in $a\in A(t,i)$.
\item[(c)] For each $(t,i)\in [0,T]\times S$, the function $c(t,i,a)$ is lower semicontinuous in $a\in A(t,i)$, and the function $\sum_{j\in S}V(j)q(j|t,i,a)$ is continuous in $a\in A(t,i)$, with  $V$ as in Condition \ref{ass:3.1}.
 \end{itemize}
\end{con}
Under Conditions \ref{ass:3.1} and  \ref{ass:3.3}(b) and (c), the function  $\sum_{j\in S}q(j|t,i,a)u(t,j)$ is continuous in $a\in A(t,i)$, for every fixed $(t,i)\in [0,T]\times S$
and $V$-bounded measurable function $u$ on $[0,T]\times S$, see the proof of Lemma 8.3.7(a) in \cite{one99}. This fact will be used in the proof of the next statement.

Also note that Condition \ref{ass:3.3} is satisfied by Example \ref{Example01}.

The main optimality result is the following one.
\begin{thm}\label{T-4.1}  Suppose Conditions \ref{ass:3.1}, \ref{ass:3.1b} and \ref{ass:3.3} are satisfied. Then there exists a unique $\varphi$ in $C_{V,V_1}^{1}([0,T]\times S)$ and some $f\in \Pi_m^d$ satisfying (\ref{E-1}) and (\ref{GLZ1}). In particular, there exists a deterministic Markov optimal policy.
\end{thm}
\par\noindent\textit{Proof.}  The statement would follow from Corollary \ref{ThC-3.2}, once we showed the existence of some $\varphi\in C_{V,V_1}^{1}([0,T]\times S)$ satisfying (\ref{E-1}). We verify this fact following a similar reasoning as in \cite{GHH2015} dealing with a risk-neutral CTMDP problem, which was also adopted in \cite{Wei:2016}, dealing with a model with a bounded cost rate. Namely, we shall obtain the desired solution $\varphi$ as a limit point of an equicontinuous family $\{\varphi_n\}$ of functions, which in turn are obtained from a sequence of CTMDP models with bounded transition and cost rates. The denumerable state space serves to prove the equicontinuity of the family $\{\varphi_n\}.$ The details are as follows.

For each integer $n\ge 1,$ let $S_n:=\{i\in S:~V(i)\le n\}$. Without loss of generality, assume for each $n\ge 1$, $S_n\ne \emptyset$.
For each $i\in S$ and $t\in [0,\infty)$, let
 $A_n(t,i):=A(t,i)$. For each $(t,i,a)\in K_{n}:=K$, define
\begin{eqnarray*}
&&q_{n}(j|t,i,a):=
q(j|t,i,a)  I_{S_{n}}(i),~\forall~j\in S,~c_n(t,i,a):=c(t,i,a)I_{S_{n}}(i),~g_n(i):=g(i)I_{S_{n}}(i).
\end{eqnarray*}
We consider the resulting sequence of CTMDP models
$
{\cal{M}}_n:=\left\{S, A_n(t,i), c_{n}, g_n, q_{n} \right\}.
$

Note that the models $\{{\cal M}_n\}$ are all with bounded transition and cost rates, and so Proposition \ref{ThCC-3.2} implies, for each $n\ge 1$, the existence of a
unique $\varphi_n$ in $C^{1}_{1,1}([0,T]\times S)$ and some $f_n\in \Pi_m^d$ satisfying
\begin{eqnarray}\label{GLZN01}
\varphi_n(s,i)-e^{g_n(i)}&=&\int_s^T \inf_{a\in A(t,i)}\left\{c_n(t,i,a)\varphi_n(t,i)+\sum_{j\in S}\varphi_n(t,j)q_n(j|t,i,a)\right\} dt \nonumber\\
&=&\int_s^T  \left\{c_n(t,i,f_n(t,i))\varphi_n(t,i)+\sum_{j\in S}\varphi_n(t,j)q_n(j|t,i,f_n(t,i))\right\}dt,\nonumber\\
&&~s\in[0,T],~i\in S.
\end{eqnarray}

Let $n\ge 1$ be fixed. For each $s\in[0,T]$, consider the $s$-shifted model
\begin{eqnarray*}
{\cal{M}}^{(s)}_n:=\left\{S, A_n^{(s)}(t,i), q_{n}^{(s)}, c_{n}^{(s)}, g_n \right\}
 \end{eqnarray*}
 with $A_n^{(s)}(t,i):=A_n(t+s,i),$ $q_{n}^{(s)}(\cdot|t,i,a):=q_{n}(\cdot|s+t,i,a)$ and $c_{n}^{(s)}(t,i,a):=c_{n}(t+s,i,a)$. Then Condition \ref{ass:3.1} is clearly satisfied by  ${\cal{M}}^{(s)}_n$, so that one can apply the reasoning in the proof of Lemma \ref{lem3.1}(c) and deduce
 \begin{eqnarray*}
 {\rm E}_i^{f_n^{(s)}}\left[e^{\int_0^{T-s} |c^{(s)}_n(t,\xi_t,f_n^{(s)}(t,\xi_t))|dt+|g_n(\xi_{T-s})|} \right]\le  Me^{T\rho}V(i)
 \end{eqnarray*}
where ${\rm E}_i^{f_n^{(s)}}$ denotes the expectation in the ${\cal{M}}^{(s)}_n$ model under the shifted policy $f_n^{(s)}(t,i):=f_n(t+s,i)$.
On the other hand, according to the uniqueness of the solution to (\ref{GLZN01}) in $C^{1}_{1,1}([0,T]\times S)$ and the discussions at the end of Section 3 of \cite{GuoZhangAMO:2017} after Theorem 3.2 therein,
 \begin{eqnarray*}
 {\rm E}_i^{f_n^{(s)}}\left[e^{\int_0^{T-s} c^{(s)}_n(t,\xi_t,f_n^{(s)}(t,\xi_t))dt+g_n(\xi_{T-s})} \right]=\varphi_n(s,i).
 \end{eqnarray*}
(The cost rate and the terminal cost were assumed to be nonnegative in \cite{GuoZhangAMO:2017}, but the results obtained there apply because ${\cal{M}}^{(s)}_n$ has bounded transition and cost rates, which can be reduced to the nonnegative case after one add to the cost rate and the terminal cost a large enough constant.) Thus, we obtain the bound
\begin{eqnarray}\label{u1}
|\varphi_n(t,i)|\leq Me^{T\rho}V(i),~\forall~ n\geq 1, (t,i)\in [0,T]\times S.
\end{eqnarray}

Next, we show that $\{\varphi_{n}, n\geq1\}$ is an equicontinuous family of functions on
$[0,T]\times S$, as follows.  Let
\begin{eqnarray*}
H_{n}(t,i):=\inf_{a\in
A_n(t,i)}\left\{c_n(t,i,a)\varphi_{n}(t,i)+\sum_{j\in S}\varphi_{n}(t,j)q_{n}(j|t,i,a)\right\},~\forall~(t,i)\in [0,T]\times S.
\end{eqnarray*}
Then, from Condition \ref{ass:3.1}  and (\ref{u1}), we see
\begin{eqnarray}\label{Li}
|H_{n}(t,i)|&\leq& \sup_{a\in A_n(t,i)}\left\{|c_n(t,i,a)\varphi_{n}(t,i)|+\sum_{j\in S}|\varphi_{n}(t,j)||q_{n}(j|t,i,a)|\right\} \nonumber \\
&\leq& \sup_{a\in A_n(t,i)}\left\{MV(i)Me^{T\rho}V(i)+Me^{T\rho}\sum_{j\in S}|q(j|t,i,a)|V(j)\right\} \nonumber \\
&\le& e^{T\rho}(M^2V^2(i)+\rho MV(i)+2M|q(i|t,i,a)|V(i)) \nonumber\\
&\leq& Me^{T\rho}M_1(3M^2+\rho)V_1(i)=:L(i),~ \forall \ (t,i)\in [0,T]\times S.
\end{eqnarray}
(Recall that $M>1$.)

Now, fix arbitrarily some $(s_0,i_0)\in [0,T]\times S$ and $\varepsilon>0$, and take
$\delta:=\min\{\frac{\varepsilon}{L(i_0)},\frac{1}{2}\}$. Then, for
every $(s,i)$ in the open neighborhood $\{(s,i)\in [0,T]\times S:~
|s-s_0|<\delta,|i-i_0|<\delta\}$, we have $i=i_0$,
and
\begin{eqnarray*}
|\varphi_{n}(s,i)-\varphi_{n}(s_{0},i_0)| &=& |\varphi_{n}(s,i_0)-\varphi_{n}(s_{0},i_0)|
=\left|\int_{s}^{T}H_{n}(t,i_0)dt-\int_{s_0}^{T}H_{n}(t,i_0)dt\right|\\
&\leq& L(i_0)|s-s_{0}|<\varepsilon,~\forall~ n\geq1.
\end{eqnarray*}
Hence, $\{\varphi_{n}, n\geq1\}$ is equicontinuous at $(s_0,i_0)$,
which, together with the arbitrariness  of $(s_0,i_0)\in [0,T]\times S$,
yields that $\{\varphi_{n}, n\geq1\}$ is equicontinuous on $[0,T]\times
S$. By Arzela-Ascoli theorem, see, e.g., p.96 of \cite{Hernandez-Lerma:1996}, there exist a subsequence $\{\varphi_{n_{k}}, k\geq1\}$ of
$\{\varphi_{n}, n\geq1\}$ and a continuous function $\varphi$ on
$[0,T]\times S$ such that
\begin{eqnarray}\label{2h}
\lim_{k\to\infty}\varphi_{n_{k}}(s,i)=\varphi(s,i),  \  {\rm and} \ |\varphi(s,i)|  \leq Me^{T\rho}V(i)~ \ \forall \ (s,i)\in [0,T]\times S,
\end{eqnarray}
where the last inequality is by (\ref{u1}).

Let \begin{eqnarray*}
H(t,i):=\inf_{a\in A(t,i)}\left\{c(t,i,a)\varphi(t,i)+\sum_{j\in
S}\varphi(t,j)q(j|t,i,a)\right\},   \forall~ (t,i)\in [0,T]\times S.
\end{eqnarray*}  We next
verify that $\lim_{k\to\infty} H_{n_k}(t,i)=H(t,i)$ for each $(t,i)\in
[0,T]\times S$, as follows. Let $(t,i)\in [0,T]\times S$ be arbitrarily fixed.  Since $q_{n_k}(j|t,i,a)\to q(j|t,i,a)$  for all $j\in S$ and $a\in A(t,i)$ as $k\to\infty$, by virtue of Lemma 8.3.7 in \cite{one99} and (\ref{u1}), we have
\begin{eqnarray*}
\limsup_{k\to \infty}H_{n_k}(t,i)&\leq& \limsup_{k\to \infty}\left\{c_{n_k}(t,i,a)\varphi_{n_k}(t,i)+\sum_{j\in S}\varphi_{n_k}(t,j)q_{n_k}(j|t,i,a)\right\} \nonumber \\
  &\leq& c(t,i,a)\varphi(t,i)+\sum_{j\in S}\varphi(t,j)q(j|t,i,a),  ~\forall~ a\in A(t,i), \nonumber
\end{eqnarray*}
so that
\begin{eqnarray}\label{H1}
\limsup_{k\to \infty}H_{n_k}(t,i)&\leq &
\inf_{a\in A(t,i)} \left\{c(t,i,a)\varphi(t,i)+\sum_{j\in
S}\varphi(t,j)q(j|t,i,a)\right\}.
\end{eqnarray}

According to the fact mentioned below Condition \ref{ass:3.3}, there exists a sequence of policies $\{f_{n_{k}}\}\subseteq
\Pi_m^d$ such that
\begin{eqnarray*}
H_{n_{k}}(t,i)&=&\inf_{a\in A(t,i)}\left\{c_{n_{k}}(t,i,a)\varphi_{n_{k}}(s,i)+\sum_{j\in S}\varphi_{n_{k}}(t,j)q_{n_{k}}(j|t,i,a)\right\} \nonumber\\
&=& c(t,i,f_{n_{k}}(t,i))\varphi_{n_{k}}(t,i)+\sum_{j\in S}\varphi_{n_{k}}(t,j)q_{n_{k}}(j|t,i,f_{n_{k}}(t,i)).
\end{eqnarray*}
Since $A(t,i)$ is
compact, by taking subsequences if necessary, we can assume without loss of generality that $\liminf_{k\rightarrow\infty}H_{n_{k}}(t,i)=\lim_{k\rightarrow\infty}H_{n_{k}}(t,i)$ and for some $a\in A(t,i)$,  $f_{n_{k}}(t,i)\to a$
as $k\to\infty$. By the virtue of Lemma 8.3.7 in \cite{one99},  we have
\begin{eqnarray*}
&&\liminf_{k\to \infty}H_{n_k}(t,i) = \liminf_{k\to \infty}\left\{ c(t,i,f_{n_{k}}(t,i))\varphi_{n_{k}}(t,i)+\sum_{j\in S}\varphi_{n_{k}}(t,j)q_{n_{k}}(j|t,i,f_{n_{k}}(t,i))\right\} \\
&\ge& c(s,i,a)\varphi(t,i)+\sum_{j\in S}\varphi(t,j)q(j|t,i,a)
\geq \inf_{a\in A(t,i)} \left\{c(t,i,a)\varphi(t,i)+\sum_{j\in S}\varphi(t,j)q(j|t,i,a)\right\}.
\end{eqnarray*}
(Recall Condition \ref{ass:3.3}.)
This, together with (\ref{H1}), implies that $\lim_{k\to\infty}
H_{n_k}(s,i)=H(s,i)$. Since $(s,i)\in[0,T]\times S$ was arbitrarily fixed, we see
from (\ref{GLZN01}), (\ref{Li}) and (\ref{2h}) that $\varphi$ satisfies
(\ref{E-1}). The same argument as in (\ref{Li}) leads to
 \begin{eqnarray*}
 |\varphi'(t,i)|=|H(t,i)|\leq   Me^{T\rho}M_1(3M^2+\rho)V_1(i),~\forall~(t,i)\in[0,T]\times S.
 \end{eqnarray*}
Therefore, we see that $\varphi\in C_{V,V_1}^{1}([0,T]\times
 S).$ The required deterministic Markov policy $f$ exists because of the fact mentioned below Condition \ref{ass:3.3}, a measurable selection theorem, see Proposition D.5 of \cite{Hernandez-Lerma:1996}.

Finally, we verify the uniqueness part. Let $\varphi\in C_{V,V_1}^{1}([0,T]\times
 S)$ be an arbitrarily fixed solution to  (\ref{E-1}). (The above reasoning shows that there exists at least one.) Let $s\in[0,T]$ be fixed, and consider the $s$-shifted model ${\cal M}^{(s)}=\left\{S,A^{(s)}(t,i), q^{(s)}, c^{(s)}, g \right\}$, which is defined as for the ${\cal M}_n^{(s)}$ model with $n$ being omitted everywhere. Let
\begin{eqnarray*}
V^{(s)}(i):=\inf_{\pi\in \Pi}{\rm E}_i^{\pi}\left[e^{\int_0^{T-s} \int_A c^{(s)}(t,\xi_t,a)\pi(da|\omega,t) dt+g(\xi_{T-s})}\right]
\end{eqnarray*}
 with ${\rm E}_i^{\pi}$ signifying the expectation in the $s$-shifted model. Then the function $\varphi^{(s)}\in C_{V,V_1}^{1}([0,T-s]\times
 S)$ defined by $\varphi^{(s)}(\tau,i):=\varphi(\tau+s,i)$ for each $(\tau,i)\in[0,T-s]\times S$ satisfies
\begin{eqnarray*}
\varphi^{(s)}(\tau,i)-e^{g(i)}&=&\int_\tau^{T-s} \inf_{a\in A^{(s)}(t,i)}\left\{c^{(s)}(t,i,a)\varphi^{(s)}(t,i)+\sum_{j\in S}\varphi^{(s)}(t,j)q^{(s)}(j|t,i,a)\right\} dt \nonumber\\
&=&\int_\tau^{T-s}  \left\{c^{(s)}(t,i,f^{(s)}(t,i))\varphi^{(s)}(t,i)+\sum_{j\in S}\varphi^{(s)}(t,j)q^{(s)}(j|t,i,f^{(s)}(t,i))\right\}dt,\nonumber\\
&&~\tau\in[0,T-s],~i\in S,
\end{eqnarray*}
for some deterministic Markov policy $f^{(s)}$. By applying Corollary \ref{ThC-3.2} to the $s$-shifted model ${\cal M}^{(s)}$, we see
$\varphi^{(s)}(0,i)=V^{(s)}(i)$, and thus $\varphi(s,i)=V^{(s)}(i)$ for each $i\in S.$ Since $s\in[0,T]$ was arbitrarily fixed, it follows that $\varphi$ is the unique solution to (\ref{E-1}) out of $\varphi\in C_{V,V_1}^{1}([0,T]\times
 S)$. The proof is completed. $\hfill\Box$
\bigskip

 \section{Conclusion}\label{sectionConclusion}
 In this paper, we considered a risk-sensitive CTMDP problem in a denumerable state space over a finite time duration. Under conditions that can be satisfied by unbounded transition and cost rates, the optimality equation was shown to have a unique solution out of a class of functions, to which Feynman-Kac formula was shown to be applicable. The results obtained in this paper can be viewed as a response to the remark in Section 7 of \cite{Wei:2016}, and complemented the relevant results in \cite{GuoZhangAMO:2017}.

\subsection*{Acknowledgement} This work is partially supported by Natural Science Foundation of Guangdong Province (Grant No.2014A030313438), Zhujiang New Star (Grant No. 201506010056), Guangdong Province outstanding young teacher training plan (Grant No. YQ2015050).

\end{document}